\documentclass[12pt]{article}
\usepackage{}
\usepackage{amsfonts}

\usepackage{amsmath, amssymb}
\usepackage{amssymb}
\def\mineappendix{
        \setcounter{section}{1}
        \setcounter{subsection}{0}
        \def\thesection{\Alph{section}}
        \def\sectionap{\@startsection  {section}{1}{\z@}
                        {-3.5ex plus-1ex minus-.2ex} {0ex plus.2ex}
                        {\reset@font\Large\bf  Appendix:  \, }
                        }
        }
\makeatother
\def\Proclaim #1. #2\par{\bigbreak\noindent{\sc#1.\enspace}{\it#2}\par}

\newcommand{\gwii}[1]{\left< \hspace{-2pt} \left< \, #1 \,
        \right>  \hspace{-2pt} \right>_{0}}

\newcommand{\gwiig}[1]{\left< \hspace{-2pt} \left< \, #1 \,
    \right> \hspace{-2pt} \right>_{g}}

\newcommand{\grava}[1]{\tau_{#1}(\gamma_{\alpha})}

\newtheorem{lem}{Lemma}[section]

\newtheorem{thm}[lem]{Theorem}

\newtheorem{rem}[lem]{Remark}

\title{ A genus-4 topological recursion relation for Gromov-Witten invariants}

\author{  Xin Wang
}

\date{}

\begin{document}
\maketitle
\begin{abstract}
In this paper, we give a new genus-4 topological recursion relation for Gromov-Witten
invariants of compact symplectic manifolds via Pixton's relations on the moduli
space of curves. As an application, we prove Pixton's relations imply a known topological
recursion relation on $\overline{\mathcal{M}}_{g,1}$ for genus $g\leq4$.
\end{abstract}
 \allowdisplaybreaks

Let $\overline{\mathcal{M}}_{g,n}$ be the moduli space of genus-$g$, stable curves with $n$ marked points. It is well known that the relations among strata only involving $\psi$ classes and boundary classes in the tautological ring of $\overline{\mathcal{M}}_{g,n}$ give universal equations for the Gromov-Witten invariants of compact symplectic manifolds. Examples of genus $g\leq3$ universal equations were given in \cite{W}, \cite{Ge1}, \cite{Ge2}, \cite{BP}, \cite{KL1} and \cite{KL2}.
Recently, a large class of relations on $\overline{\mathcal{M}}_{g,n}$, called Pixton's relations, was proposed in \cite{P1} and proved in \cite{PPZ}. They are a kind of combinatorial formulas among strata involving $\kappa$ classes, $\psi$ classes and boundary classes. Pixton wrote a Sage program in \cite{P2} to compute these finitely many relations for fixed $g$ and $n$. We will use it to eliminate all the $\kappa$ classes appearing in Pixton's formula and obtain the kind of relations we want.
 In this paper, we give a new genus-4 universal equation, i.e. a genus-4 topological recursion relation on $\overline{\mathcal{M}}_{4,1}$.

 To describe universal equations, we need some conventions. Let $M$ be any compact symplectic manifold. Define the {\it big phase space}  for Gromov-Witten invariants of $M$ to be the product of infinitely many copies of $H^{*}(M;\mathbb{C})$. We will fix a basis $\{\gamma_{\alpha}|\alpha=1,...,N\}$ of $H^{*}(M;\mathbb{C})$.
 An operator $T$ on the space of vector fields on the big phase space was introduced
 in \cite{L1} to simplify topological recursion relations. The operator $T$ is very useful when we translate relations in the tautological ring of $\overline{\mathcal{M}}_{g,n}$ into universal  equations for Gromov-Witten invariants. We will write the universal equations of Gromov-Witten invariants as equations among correlators $\langle\hspace{-2pt}\langle \mathcal{W}_{1}\cdot\cdot\cdot \mathcal{W}_{k}\rangle\hspace{-2pt}\rangle_{g}$ which are by definition the $k$-th covariant derivatives of the generating functions of genus-$g$ Gromov-Witten invariants with respect to the trivial connection on the big phase space. we will briefly review these definitions in Section 1 for completeness.

 In this paper, we will prove the following genus-4 universal equation
 \begin{thm}
 \label{TRR4}
 For Gromov-Witten invariants of any compact symplectic manifold, the following topological recursion relation holds for any vector field $\mathcal{W}$ on the big phase space:
  \begin{align}
 \langle\hspace{-2pt}\langle T^{4}(\mathcal{W}) \rangle\hspace{-2pt}\rangle _{4}=A(\mathcal{W})
 \label{genus4TRR}
\end{align}
where $A(\mathcal{W})$ is given by a very complicated  formula  only involving invariants of genus not bigger than 3.
The explicit expression of $A(\mathcal{W})$ can be found in Appendix and we will not put it here.
\end{thm}
Equivalently, this theorem corresponds to a 
relation in the tautological ring of $\overline{\mathcal{M}}_{4,1}$, representing $\psi_{1}^{4}$ as a linear combination of boundary strata decorated with $\psi$ classes. Since the relation is very long, we will not show its dual graph representations in this paper. The interested readers  can write it carefully by themselves. Actually, our main idea behind the proof of Theorem 0.1 is: first obtain the relation representing $\psi_{1}^{4}$ as a linear combination of strata containing $\psi$ classes and boundary classes on $\overline{\mathcal{M}}_{4,1}$ via Pixton's relation, then translate it to equation \eqref{genus4TRR} by the splitting principle of Gromov-Witten invariants.

We should mention that in \cite{L2} it was conjectured that the following type of topological recursion relations hold for all genera $g$,
  \begin{align}
\left<\hspace{-2pt}\left<T^{g}(\mathcal{W})\right>\hspace{-2pt}\right>_{g}=A_{g-1}(\mathcal{W}).
\label{all g}
\end{align}
Equivalently, equation \eqref{all g} gives a topological recursion relation on $RH^{2g}(\overline{\mathcal{M}}_{g,1})$, representing $\psi_{1}^{g}$ as a linear combination of boundary strata with $\psi$ classes.
It is well known that \eqref{all g} holds for $g\leq3$ (cf. \cite{W}, \cite{Ge2} and \cite{KL1}). Thus our Theorem 0.1 gives an affirmative answer to this conjecture at genus-4. In general, if it exists,
such kind of relations
goes very complicated as the genus $g$ becomes more and more complicated. However,  very interesting topological recursion relations on  $RH^{2(2g+r)}(\overline{\mathcal{M}}_{g,1})$ (for any $r\geq0$) were found and proved in \cite{LP} for all genus $g\geq1$ by localisation techniques.
\begin{align}
\psi_{1}^{2g+r}=\sum_{g_{1}+g_{2}=g,g_{i}>0}\sum_{a+b=2g-1+r}(-1)^{a}\frac{g_{2}}{g}\rho_{*}
\Big(\psi_{*1}^{a}\psi^{b}_{*2}\cap[\Delta_{1,\emptyset}(g_{1},g_{2})]\Big).
\label{g,1}
\end{align}
Given that Pixton's relations
are conjectured to generate all relations in the tautological ring,
one would expect they imply the relations \eqref{g,1}. We prove this in the case of genus
 $g \leq4$.
\begin{thm}
In the case of genus
 $g \leq4$, Pixton's relations imply  the topological recursion relation \eqref{g,1}.
\end{thm}
As for genus $g\geq5$, we conjecture that  our method in the proof also works as long as we get the topological recursion relation on $RH^{2g}(\overline{\mathcal{M}}_{g,1})$.

This paper is organized as follows.  In Section 1, we introduce notations and review
basic theories needed in this paper. We prove Theorem 0.1 in Section 2 and prove Theorem 0.2
 in Section 3. In the Appendix, we give the
precise formula of the topological recursion relation \eqref{genus4TRR}.
\\{\bf Acknowledgements.}
The author is partially supported by  the fundamental Research Funds of Shandong
University. We would like to thank professor Xiaobo Liu for helpful discussion. We  also thank Felix Janda for lots of helpful comments on the paper.

\section{Preliminaries}

\subsection{Recollection of  Pixton's relations on moduli space of curves}

Let $\overline{\mathcal {M}}_{g,n}$ be the moduli space of stable curves of genus-$g$ with $n$-marked points. The boundary strata of the moduli $\overline{\mathcal {M}}_{g,n}$ of fixed topological type correspond to stable graphs.

 Formally, a stable graph is the structure
 $$\Gamma=(V, E, L ,g)$$
 satisfying the following properties:
 \begin{itemize}
\item{$V$ is the vertex set with a genus function $g: V\rightarrow \mathbb{Z}$,}
\item{$E$ is the edge set,}
\item{$L$ is the set of legs (corresponding  to the set of markings),}
\item{the pair $(V,E)$ defines a connected graph,}
\item{for each vertex $v$, the stability condition holds:
$$2g(v)-2+n(v)>0,$$ where $n(v)$ is the degree of the vertex $v$ in the graph $\Gamma$.
}
\end{itemize}

 The genus of a stable graph is defined by
 $$g(\Gamma):=\sum_{v\in V}g(v)+h^{1}(\Gamma).$$
 To each stable graph $\Gamma$, we associate the moduli space
 $$\overline{\mathcal {M}}_{\Gamma}:=\prod_{v\in V}\overline{\mathcal {M}}_{g(v),n(v)}.$$
 Then there is a canonical morphism
 $$\xi_{\Gamma}: \overline{\mathcal {M}}_{\Gamma}\rightarrow \overline{\mathcal {M}}_{g,n}$$

Recall that the {\it strata algebra} $\mathcal{S}_{g,n}$ on $\overline{\mathcal {M}}_{g,n}$ is defined to be a finite dimensional algebra, with additive basis the isomorphism classes $[\Gamma,\gamma]$, where $\Gamma$ is a stable graph of genus $g$ with $n$ legs and $\gamma$ is a class on   $\overline{\mathcal {M}}_{\Gamma}$ which is a product of monomials in $\kappa$ classes on each vertex of $\Gamma$ and $\psi$ classes on each half edge.
Push forward along $\xi_{\Gamma}$ defines a  ring homomorphism from strata algebra to cohomology ring
\begin{align*}
q: &\mathcal{S}_{g,n}\rightarrow H^{*}(\overline{\mathcal {M}}_{g,n})
\\& [\Gamma, \gamma]\rightarrow (\xi_{\Gamma})_{*}(\gamma).
\end{align*}
The {\it tautological ring} $RH^{*}(\overline{\mathcal {M}}_{g,n})$ is by definition
the image of $q$. An element in the kernel of $q$ is call a {\it tautological relation}
on $\overline{\mathcal {M}}_{g,n}$.

We define tautological classes $\mathcal {R}(g,n,r;\sigma, a_{1},...,a_{n})$ associated to the data:
   \begin{itemize}
\item{$g,n\in \mathbb{Z}_{\geq0}$ in the stable range $2g-2+n>0$,}
\item{$a_{1},...,a_{n}$ are nonnegative integers not 2 mod 3,}
\item{$\sigma$: a partition with no parts 2 mod 3,}
\item{$r\in \mathbb{Z}_{\geq0}$ satisfying
$$|\sigma|+\sum_{i}a_{i}\leq 3r-g-1,$$
$$|\sigma|+\sum_{i}a_{i}\equiv 3r-g-1 \pmod{2} $$}
where partition $\sigma=1^{n_{1}}3^{n_{3}}4^{n_{4}}6^{n_{6}}\cdot\cdot\cdot$, $|\sigma|:=\sum_{i\neq2 mod3}in_{i}$.
\end{itemize}
The elements $\mathcal {R}(g,n,r;\sigma, a_{1},...,a_{n})$ are expressed as sums over stable graphs of genus $g$ with $n$ legs.
\begin{align}
 \mathcal {R}(g,n,r;\sigma, a_{1},...,a_{n}):=
 \sum_{\Gamma\in G_{g,n}}\frac{1}{|Aut(\Gamma)|} (\xi_{\Gamma})_{*}(\mathcal {R}_{\Gamma}(g,n,r;\sigma, a_{1},...,a_{n}))
 \label{sumgra}
 \end{align}
 where $G_{g,n}$ denotes the finite set of stable graphs of genus $g$ with $n$ legs (up to isomorphism).
Pixton's relations then takes 
\begin{align}
\mathcal {R}(g,n,r;\sigma, a_{1},...,a_{n})=0 \in H^{2r}(\overline{\mathcal {M}}_{g,n}, \mathbb{Q}).
\label{Pixrel}
\end{align}

Before writing the formula for $\mathcal {R}_{\Gamma}(g,n,r;\sigma, a_{1},...,a_{n})$, a few definitions are required.
We first introduce the following two power series:
\begin{align*}
A(T)=\sum_{n\geq0}\frac{(6n)!}{(3n)!(2n)!}T^{n},
\end{align*}
\begin{align*}
B(T)=\sum_{n\geq0}\frac{6n+1}{6n-1}\frac{(6n)!}{(3n)!(2n)!}T^{n}.
\end{align*}
These two series play a central role in Pixton's relations. Let
\begin{align*}
\widehat{C}_{3i}(T,\zeta)=T^{i}A(\zeta T)
\end{align*}
\begin{align*}
\widehat{C}_{3i+1}(T,\zeta)=\zeta T^{i}B(\zeta T)
\end{align*} where $\zeta$ is a formal parameter satisfying $\zeta^{2}=1$.

For any polynomial in variable $T$,
\begin{align*}
F(T)=\sum_{n}a_{n}T^{n}+\zeta\sum_{n}b_{n}T^{n},
\end{align*}
define
\begin{align*}
\{F(T)\}:=\sum_{n}a_{n}T^{n}K_{n,0}+\sum_{n}b_{n}T^{n}K_{n,1}.
\end{align*}
For each vertex $v\in V$ of a stable graph, we introduce an auxiliary variable $\zeta_{v}$  and impose the conditions
$$\zeta_{v}\zeta_{v'}=\zeta_{v'}\zeta_{v},\quad   \zeta_{v}^{2}=1.$$
The variables $\zeta_{v}$ will be responsible for keeping track of a local parity condition at each vertex.

For each stable graph $\Gamma\in G_{g,n}$,
   define
\begin{align}
\widehat{\kappa}_{\Gamma}(K_{e_{1},a_{1}}\cdot\cdot\cdot K_{e_{l},a_{l}}):=
\sum_{\tau\in S_{l}}\prod_{c \text{cycle in} \tau} (\sum_{v\in V(\Gamma)}
\kappa_{\sum_{i\in c}e_{i}}^{(v)}\zeta_{v}^{\sum_{i\in c}a_{i}})
\end{align}

%
%
%

We denote by $\mathcal {R}_{\Gamma}(g,n,r;\sigma, a_{1},...,a_{n}) \in H^{2r}(\overline{\mathcal {M}}_{g,n}; \mathbb{Q}) $ the class
\begin{align*}
&\mathcal {R}_{\Gamma}(g,n,r;\sigma,a_{1},...,a_{n})
\nonumber\\=&\frac{1}{2^{h^{1}(\Gamma)}}
\Big[
\widehat{\kappa}_{\Gamma}\big(e^{\{1-\widehat{C}_{0}\}}\{\widehat{C}_{\sigma_{1}}\}
\cdot\cdot\cdot\{\widehat{C}_{\sigma_{l}}\}\big)
\prod_{i=1}^{n}\widehat{C}_{a_{i}}(\psi_{i}T,\zeta_{v_{i}})
\prod_{e\in E}\Delta_{e}\Big]_{T^{r-|E|}\prod_{v}\zeta_{v}^{g(v)+1}}
\end{align*}
where $h^{1}(\Gamma)=|E(\Gamma)|-|V(\Gamma)|+1$ is the number of loops in $\Gamma$, marking $i$ corresponds to half-edge $h_{i}$ on vertex $v_{i}$, and $g(v)$ is the genus of vertex $v$.
The subscript $T^{r-|E|}\prod_{v}\zeta_{v}^{g(v)+1}$ indicates the coefficient of the monomial $T^{r-|E|}\prod_{v}\zeta_{v}^{g(v)+1}$ after the product inside the brackets is expanded. For each edge $e\in E$,
    \begin{align*}
\Delta_{e}=
\frac{A(\zeta_{1}\psi_{1}T)\zeta_{2}B(\zeta_{2}\psi_{2}T)+
\zeta_{1}B(\zeta_{1}\psi_{1}T)A(\zeta_{2}\psi_{2} T)+\zeta_{1}+\zeta_{2}}{(\psi_{1}+\psi_{2})T},
\end{align*}
where $\zeta_{1},\zeta_{2}$ are the $\zeta$ variables assigned to the vertices adjacent to the edge $e$ and $\psi_{1}, \psi_{2}$ are the $\psi$ classes corresponding to the half edges.
    The numerator of $\Delta_{e}$  is divisible by the denominator due to the identity
    \begin{align*}
A(T)B(-T)+A(-T)B(T)=-2.
\end{align*}


\subsection{Basics on Gromov-Witten invariants}
Let $M$ be a compact symplectic manifold. The {\it small phase space} is by
definition $H^{*}(M;\mathbb{C})$ and the {\it big phase space} is defined to be
$\mathcal{P}:=\prod_{n=0}^{\infty} H^{*}(M,\mathbb{C})$.
Let $\{\gamma_{1},...,\gamma_{N}\}$ be a fixed basis of $H^{*}(M,\mathbb{C})$, where $\gamma_{1}$
is the identity element of the cohomology ring of $M$.
The  corresponding
basis for the $n$-th copy of $H^{*}(M,\mathbb{C})$ in this product is denoted by
$\{ \grava{n} \mid \alpha=1, \ldots, N\}$ for $ n \geq 0$.
Let $t_{n}^{\alpha}$ be the coordinates on $\mathcal{P}$
with respect to the standard basis $\{ \grava{n} \mid \alpha=1, \ldots, N, \,\,\, n \geq 0\}$.
%
%
%
We will identify $\tau_{n}(\gamma_{\alpha})$ with $\frac{\partial}{\partial t_{n}^{\alpha}}$ as vector fields on the big phase space. If $n<0$, $\tau_{n}(\gamma_{\alpha})$ is understood to be the zero vector field. We will also write $\tau_{0}(\gamma_{\alpha})$ simply as $\gamma_{\alpha}$.
We use $\tau_{+}$ and $\tau_{-}$ to denote the operators which shift the level of descendants, i.e.
$$\tau_{\pm}\Big(\sum_{n,\alpha}f_{n,\alpha}\tau_{n}(\gamma_{\alpha})\Big)
=\sum_{n,\alpha}f_{n,\alpha}\tau_{n\pm1}(\gamma_{\alpha})$$ where $f_{n,\alpha}$ are functions on the big phase space.

 Define $\eta=(\eta_{\alpha\beta})$ to be the matrix of intersection pairing on $H^{*}(M;\mathbb{C})$  in this basis $\{\gamma_{1},...,\gamma_{N}\}$. We will use
$\eta=(\eta_{\alpha\beta})$ and $\eta^{-1}=(\eta^{\alpha\beta})$ to lower and  raise indices, for example $\gamma^{\alpha}:=\eta^{\alpha\beta}\gamma_{\beta}$ for any $\alpha$.
Here we use the summation convention that repeated indices should be summed over their entire ranges.

Let
$$\left<\tau_{n_{1}}(\gamma_{\alpha_{1}})\cdot\cdot\cdot \tau_{n_{k}}(\gamma_{\alpha_{k}})\right>_{g,\beta}
:=\int_{[\overline{\mathcal{M}}_{g,n}(M,\beta)]^{vir}}\prod_{i=1}^{k}(\Psi_{i}\cup ev_{i}^{*}(\gamma_{\alpha_{i}})) $$
be the genus-$g$, degree-$\beta$, descendant Gromov-Witten invariant associated to $\gamma_{\alpha_{1}},...,\gamma_{\alpha_{k}}$ and nonnegative integers $n_{1},...,n_{k}$ (cf. \cite{RT}, \cite{LT}). Here $\overline{\mathcal{M}}_{g,n}(M,\beta)$ is the moduli space of stable maps from genus-$g$, $k$-marked curves to $M$ of degree $\beta\in H_{2}(M;\mathbb{Z})$. $\Psi_{i}$ is the first Chern class of the tautological line bundle over $\overline{\mathcal{M}}_{g,n}(M,\beta)$ whose geometric fiber is the cotangent space of the domain curve at the $i$-th marked point and $ev_{i}: \overline{\mathcal{M}}_{g,n}(M,\beta)\rightarrow M$ is the $i$-th evaluation map for all $i=1,...,k$ and $[\overline{\mathcal{M}}_{g,n}(M,\beta)]^{vir}$ is the virtual fundamental class. The genus-$g$ generating function is defined to be
\begin{align*}
F_{g}:=\sum_{k\geq0}\sum_{\alpha_{1},...,\alpha_{k}}
\sum_{n_{1},...,n_{k}}
\frac{1}{k!}t_{n_{1}}^{\alpha_{1}}\cdot\cdot\cdot t_{n_{k}}^{\alpha_{k}}\sum_{\beta}q^{\beta}
\left<\tau_{n_{1}}(\gamma_{\alpha_{1}})\cdot\cdot\cdot \tau_{n_{k}}(\gamma_{\alpha_{k}})\right>_{g,\beta}
\end{align*}
where $q^{\beta}$ belongs to the Novikov ring. This function is understood as a formal power series in variables $\{t_{n}^{\alpha}\}$ with coefficients in the Novikov ring.

Define a $k$-tensor $\gwiig{\cdot\cdot\cdot}$ by
\begin{align}
\gwiig{\mathcal{W}_{1}\mathcal{W}_{2}\cdot\cdot\cdot\mathcal{W}_{k}}
:=\sum_{m_{1},\alpha_{1},...,m_{k},\alpha_{k}}f_{m_{1},\alpha_{1}}^{1}
\cdot\cdot\cdot f_{m_{k},\alpha_{k}}^{k}\frac{\partial^{k}}{\partial t_{m_{1}}^{\alpha_{1}}\cdot\cdot\cdot \partial t_{m_{k}}^{\alpha_{k}}}F_{g}
\end{align}
for vector fields $\mathcal{W}_{i}=\sum_{m,\alpha}f_{m,\alpha}^{i}\frac{\partial}{\partial t_{m}^{\alpha}}$ where $f_{m,\alpha}^{i}$ are functions on the big phase space. This tensor is called $k$-$point$ $(correlation)$ $function$.

For any vector fields $\mathcal{W}_{1}$ and $\mathcal{W}_{2}$ on the big phase space, the $quantum$ $product$ of $\mathcal{W}_{1}$ and $\mathcal{W}_{2}$ is defined by
$$\mathcal{W}_{1}\circ \mathcal{W}_{2}:=\gwii{\mathcal{W}_{1}\mathcal{W}_{2}
\gamma^{\alpha}}\gamma_{\alpha}.$$
Define the vector field
$$T(\mathcal{W}):=\tau_{+}(\mathcal{W})-\langle\hspace{-2pt}\langle \mathcal{W}\gamma^{\alpha}\rangle\hspace{-2pt}\rangle_{0}\gamma_{\alpha}$$
for any vector field $\mathcal{W}$. The operator $T$ was introduced in \cite{L1} in order to simplify topological recursion relations of Gromov-Witten invariants. Let $\nabla$ be the trivial flat connection on the big phase space with respect to the coordinates $\{t_{n}^{\alpha}\}$. Then the covariant derivative of the quantum product is given by
\begin{align}
\nabla_{\mathcal{W}_{3}}(\mathcal{W}_{1}\circ \mathcal{W}_{2})=(\nabla_{\mathcal{W}_{3}}\mathcal{W}_{1})\circ \mathcal{W}_{2}+\mathcal{W}_{1}\circ(\nabla_{\mathcal{W}_{3}} \mathcal{W}_{2})+\langle\hspace{-2pt}\langle \mathcal{W}_{1}\mathcal{W}_{2}\mathcal{W}_{3}\gamma^{\alpha}\rangle\hspace{-2pt}\rangle_{0} \gamma_{\alpha}
\label{derquan}
\end{align}
and the covariant derivative of the operator $T$ satisfies
\begin{align}
\nabla_{\mathcal{W}_{2}}T(\mathcal{W}_{1})=T(\nabla_{\mathcal{W}_{2}}\mathcal{W}_{1})-\mathcal{W}_{1}\circ \mathcal{W}_{2}
\label{derT}
\end{align}
for any vector fields $\mathcal{W}_{1}$, $\mathcal{W}_{2}$ and $\mathcal{W}_{3}$. We need these formulas to compute derivatives of universal equations.
%

\section{Proof of Theorem 0.1 }
Obviously, Pixton's relations give  lots of tautological relations on $\overline{\mathcal {M}}_{g,n}$ for any $g$ and $n$. What we need is to  eliminate the kappa classes in Pixton's  relations to get explicit topological recursion relations on $\overline{\mathcal {M}}_{g,n}$.

The proof of Theorem 0.1 mainly contains five steps. Here we should mention that, by using his sage code, Pixton has checked the presence of Getzler's relation \cite{Ge1} in $RH^{4}(\overline{\mathcal {M}}_{1,4})$ and   Belorousski-Pandharipande relation \cite{BP} in $RH^{4}(\overline{\mathcal {M}}_{2,3})$ (cf. \cite{P1}). 

Actually, the
equation \eqref{genus4TRR} is also obtained with the aid of Pixton's sage code. In principle, we can use this method for all $g$ and $n$ to get topological recursion relations. Here we focus on relations on $\overline{\mathcal {M}}_{4,1}$.

Firstly, find all the stable graphs up to isomorphism. In fact, we can define a degenerate operation on any vertex of a stable graph: simply splitting one vertex with genus $g(v)$ into two vertexes with genus  $g(v_{1})$ and  $g(v_{2})$ respectively,  such that $g(v)=g(v_{1})+g(v_{2})$.  We notice that any stable graph $\Gamma$ of genus-$g=\sum_{v\in V}g(v)+h^{1}(\Gamma)$ with $n$-legs can be obtained from a graph with only one vertex attached with genus-$\sum_{v\in V(\Gamma)}g(v)$, $h^{1}(\Gamma)$ loops and $n$-legs via finite (equal to the number of ordinary edges of $\Gamma$) steps of degeneration.
Then we identify graphs which are isomorphic to each other.

Secondly,  we decorate each stable graph with some $\kappa$ classes on the vertex and $\psi$ classes on the half edge and legs, and identify the isomorphic ones. So we get the canonical linear basis of strata algebra $\mathcal{S}_{g,n}$.
Then we reorder the basis so that strata with  $\kappa$ classes are in the front.

Thirdly, let $R_{g,n}$ be the set of all elements in the strata algebra $\mathcal{S}_{g,n}$ produced as follows: choose a dual graph $\Gamma$  for a boundary stratum of $\overline{\mathcal{M}}_{g,n}$, pick one of the components $\mathcal{S}_{g',n'}$ in
$\mathcal{S}_{\Gamma}=\mathcal{S}_{g_{1},n_{1}}\times...\times\mathcal{S}_{g_{m},n_{m}}$, take the product of a relation $\mathcal{R}(g',n',d;\sigma,a_{1},...,a_{n})=0$ on the chosen component
together with arbitrary classes on the other components, and push forward along the gluing map
$\mathcal{S}_{\Gamma}\rightarrow \mathcal{S}_{g,n}$. So the set $R_{g,n}$ give us linear equations between the generators of the strata algebra. For our purpose,  we should rewrite all these linear equations  with respect to our list of strata and get the corresponding coefficients matrix $R$.

Fourthly, by technics in linear algebra, we can transform the matrix $R$ into the reduced row  echelon form. Then we can get many linear independent tautological relations among strata only involving $\psi$ classes and boundary classes.

Lastly, we translate the tautological relations obtained above into differential equations for generating functions of Gromov-Witten invariants, here each $\psi$ class corresponds to an insertion of the operator $T$. Actually we can get many linear independent differential equations, from which we can choose the  equation \eqref{genus4TRR}.
\begin{rem}
Actually, the first 3 steps in the above algorithm is essentially the same as the method in Pixton's sage program except the reordering of the strata generators.
Using this method, we have also checked for the presence of all known topological recursion relations (cf. \cite{Ge2}, \cite{KL1} and  \cite{KL2}). To our knowledge, Lin and Zhou have done some similar work (cf. \cite{Li}). On the other hand, we get many new relations for genus-2 and genus-3 and  their relation with higher genus Virasoro conjecture will be studied  in the forthcoming papers.
\end{rem}
\section{Proof of Theorem 0.2 }
First we use operator $T$ to translate equation \eqref{g,1} into a differential equation
 \begin{align}
 \langle\hspace{-2pt}\langle T^{2g+r}(\mathcal{W})\rangle\hspace{-2pt}\rangle_{g}
 =\sum_{g_{1}+g_{2}=g,g_{i}>0}\sum_{a+b=2g-1+r}(-1)^{a}\frac{g_{2}}{g}
  \langle\hspace{-2pt}\langle \mathcal{W} T^{a}(\gamma^{\alpha})\rangle\hspace{-2pt}\rangle_{g_{1}}
   \langle\hspace{-2pt}\langle T^{b}(\gamma_{\alpha})\rangle\hspace{-2pt}\rangle_{g_{2}}
   \label{eqXPrel}
 \end{align}
 where $g\geq1$, $r\geq0$ and $\mathcal{W}$ is an arbitrary vector field on the big phase space.

 Via the operator $T$, we can reformulate the genus-0 topological recursion relation as
 \begin{align}
 \langle\hspace{-2pt}\langle T(\mathcal{W}_{1})\mathcal{W}_{2}\mathcal{W}_{3}\rangle\hspace{-2pt}\rangle_{0}=0
 \label{genus0TRR}
 \end{align}
 and the genus-1 topological recursion relation as
 \begin{align}
 \langle\hspace{-2pt}\langle T(\mathcal{W})\rangle\hspace{-2pt}\rangle_{1}=\frac{1}{24}\langle\hspace{-2pt}\langle \mathcal{W}\gamma^{\alpha}\gamma_{\alpha}\rangle\hspace{-2pt}\rangle_{0}
 \label{genus1TRR}
 \end{align}  for any vector fields $\mathcal{W}$ and $\mathcal{W}_{i}$.

 Taking derivatives of equation \eqref{genus1TRR} and combining with equation \eqref{derT}, we have
  \begin{align*}
 \langle\hspace{-2pt}\langle T(\mathcal{W})\mathcal{V}\rangle\hspace{-2pt}\rangle_{1}=\langle\hspace{-2pt}\langle \{\mathcal{W}\circ \mathcal{V} \}\rangle\hspace{-2pt}\rangle_{1}+\frac{1}{24}\langle\hspace{-2pt}\langle \mathcal{W}\mathcal{V}\gamma^{\alpha}\gamma_{\alpha}\rangle\hspace{-2pt}\rangle_{0}.
 \end{align*}
 To prove theorem 0.2, we only need to consider the following nontrivial cases.
 \subsection{Case of genus $g=2$}
 Due to the dimension constraints, we only need to show that equation~\eqref{eqXPrel} holds for $r=0$, i.e.
  \begin{align}
 \langle\hspace{-2pt}\langle T^{4}(\mathcal{W})\rangle\hspace{-2pt}\rangle_{2}
 =\frac{1}{2}
  \langle\hspace{-2pt}\langle \mathcal{W} T^{2}(\gamma^{\alpha})\rangle\hspace{-2pt}\rangle_{1}
   \langle\hspace{-2pt}\langle T(\gamma_{\alpha})\rangle\hspace{-2pt}\rangle_{1}.
 \label{g2LP}
 \end{align}

 We recall the genus-2 Mumford relation (cf. \cite{Ge2}) as formulated in \cite{L1}:
 \begin{align}
 &\langle\hspace{-2pt}\langle T^{2}(\mathcal{W})\rangle\hspace{-2pt}\rangle_{2}\nonumber
\\=&\frac{7}{10}{\langle\hspace{-2pt} \langle {\gamma ^\alpha }\rangle \hspace{-2pt}\rangle _1}{\langle \hspace{-2pt}\langle \{\gamma_\alpha \circ \mathcal{W}\}\rangle \hspace{-2pt}\rangle _1}+\frac{1}{10}{\langle \hspace{-2pt}\langle {\gamma ^\alpha }\{\gamma_\alpha \circ \mathcal{W}\}\rangle \hspace{-2pt}\rangle _1}-\frac{1}{240}{\langle \hspace{-2pt}\langle \mathcal{W}\{{\gamma ^\alpha }\circ \gamma_\alpha\}\rangle \hspace{-2pt}\rangle _1}\nonumber
\\&+\frac{13}{240}{\langle \hspace{-2pt}\langle \mathcal{W}{\gamma^\alpha } \gamma_\alpha \gamma^{\beta}\rangle\hspace{-2pt} \rangle _0}
\langle\hspace{-2pt}\langle\gamma_{\beta}\rangle\hspace{-2pt}\rangle_{1}
+\frac{1}{960}{\langle\hspace{-2pt} \langle \mathcal{W}{\gamma^\alpha } \gamma_\alpha \gamma^{\beta}\gamma_{\beta}\rangle \hspace{-2pt}\rangle _0}
\label{genus2TRR}
 \end{align}
 for any vector field $\mathcal{W}$.

By equation~\eqref{genus2TRR}, \eqref{genus0TRR} and \eqref{genus1TRR} and  their  derivatives, both sides of equation~\eqref{g2LP} equal to the expression
    \begin{align*}
\frac{1}{1152}
  \langle\hspace{-2pt}\langle \Delta\Delta \mathcal{W} \rangle\hspace{-2pt}\rangle_{0}.
\end{align*}
\subsection{Case of genus $g=3$}
  Due to the dimension constraints, we only need to equation \eqref{eqXPrel} holds for $r=0$ and $r=1$. For simplicity, we here only give the proof of the case for $r=0$ and the other case can be similarly done, i.e.
  \begin{align}
 \langle\hspace{-2pt}\langle T^{6}(\mathcal{W})\rangle\hspace{-2pt}\rangle_{3}
   =\frac{1}{3}{\langle\hspace{-2pt}\langle \mathcal{W}{T^4}({\gamma ^\alpha })\rangle\hspace{-2pt}\rangle _2}{\langle\hspace{-2pt}\langle T({\gamma _\alpha })\rangle\hspace{-2pt}\rangle _1} - \frac{1}{3}{\langle\hspace{-2pt}\langle \mathcal{W}{T^5}({\gamma ^\alpha })\rangle\hspace{-2pt}\rangle _2}{\langle\hspace{-2pt}\langle {\gamma _\alpha }\rangle\hspace{-2pt}\rangle _1}.
   \label{g3LP}
 \end{align}
 To prove equation~\eqref{g3LP}, we need the genus-3 topological recursion relation (cf. \cite{KL1}):
 \begin{align}
&\langle\hspace{-2pt}\langle T^{3}(\mathcal{W})\rangle\hspace{-2pt}\rangle_{3}\nonumber
\\=&-\frac{1}{252}\langle\hspace{-2pt}\langle \mathcal{W}T(\gamma_{\alpha}\circ \gamma^{\alpha})\rangle\hspace{-2pt}\rangle_{2}+\frac{5}{42}\langle\hspace{-2pt}\langle T(\gamma_{\alpha})\{\mathcal{W}\circ \gamma^{\alpha}\}\rangle\hspace{-2pt}\rangle_{2}\nonumber
\\&+\frac{13}{168}\langle\hspace{-2pt}\langle T(\gamma^{\alpha})\rangle\hspace{-2pt}\rangle_{2}\langle\hspace{-2pt}\langle\gamma_{\alpha}\mathcal{W}\gamma^{\beta}\gamma_{\beta}\rangle\hspace{-2pt}\rangle_{0}
+\frac{41}{21}\langle\hspace{-2pt}\langle T(\gamma^{\alpha})\rangle\hspace{-2pt}\rangle_{2}\langle\hspace{-2pt}\langle\{\gamma_{\alpha}\circ \mathcal{W}\}\rangle\hspace{-2pt}\rangle_{1}
\nonumber
\\&-\frac{13}{168}\langle\hspace{-2pt}\langle\{\mathcal{W}\circ \gamma_{\alpha}\circ \gamma^{\alpha}\}\rangle\hspace{-2pt}\rangle_{2}
+\frac{1}{280}\langle\hspace{-2pt}\langle \mathcal{W}\gamma^{\alpha}\rangle\hspace{-2pt}\rangle_{1}\langle\hspace{-2pt}\langle\gamma_{\alpha}\{\gamma^{\beta}\circ \gamma_{\beta}\}\rangle\hspace{-2pt}\rangle_{1}
\nonumber
\\&-\frac{23}{5040}\langle\hspace{-2pt}\langle\gamma^{\alpha}\rangle\hspace{-2pt}\rangle_{1}\langle\hspace{-2pt}\langle\gamma_{\alpha}\mathcal{W}\{\gamma^{\beta}\circ\gamma_{\beta}\}\rangle\hspace{-2pt}\rangle_{1}
-\frac{47}{5040}\langle\hspace{-2pt}\langle\gamma^{\alpha}\rangle\hspace{-2pt}\rangle_{1}\langle\hspace{-2pt}\langle\gamma_{\alpha}\gamma^{\beta}\rangle\hspace{-2pt}\rangle_{1}
\langle\hspace{-2pt}\langle\gamma_{\beta}\mathcal{W}\gamma^{\mu}\gamma_{\mu}\rangle\hspace{-2pt}\rangle_{0}
\nonumber
\\&-\frac{5}{1008}\langle\hspace{-2pt}\langle \mathcal{W}\gamma^{\alpha}\rangle\hspace{-2pt}\rangle_{1}
\langle\hspace{-2pt}\langle\gamma_{\alpha}\gamma^{\beta}\gamma_{\beta}\gamma^{\mu}\rangle\hspace{-2pt}\rangle_{0}\langle\hspace{-2pt}\langle\gamma_{\mu}\rangle\hspace{-2pt}\rangle_{1}
+\frac{23}{504}\langle\hspace{-2pt}\langle\gamma^{\alpha}\rangle\hspace{-2pt}\rangle_{1}\langle\hspace{-2pt}\langle\gamma_{\alpha}\mathcal{W}\gamma^{\beta}\gamma_{\beta}\gamma^{\mu}\rangle\hspace{-2pt}\rangle_{0}
\langle\hspace{-2pt}\langle\gamma_{\mu}\rangle\hspace{-2pt}\rangle_{1}
\nonumber
\\&+\frac{11}{140}\langle\hspace{-2pt}\langle\gamma^{\alpha}\gamma^{\beta}\rangle\hspace{-2pt}\rangle_{1}\langle\hspace{-2pt}\langle\gamma_{\alpha}\{\gamma_{\beta}\circ \mathcal{W}\}\rangle\hspace{-2pt}\rangle_{1}-\frac{4}{35}\langle\hspace{-2pt}\langle\gamma^{\alpha}\rangle\hspace{-2pt}\rangle_{1}\langle\hspace{-2pt}\langle\gamma_{\alpha}\gamma^{\beta}\rangle\hspace{-2pt}\rangle_{1}\langle\hspace{-2pt}\langle\{\gamma_{\beta}\circ \mathcal{W}\}\rangle\hspace{-2pt}\rangle_{1}
\nonumber
\\&+\frac{2}{105}\langle\hspace{-2pt}\langle \mathcal{W}\gamma^{\alpha}\rangle\hspace{-2pt}\rangle_{1}\langle\hspace{-2pt}\langle\{\gamma_{\alpha}\circ\gamma_{\beta}\}\rangle\hspace{-2pt}\rangle_{1}\langle\hspace{-2pt}\langle\gamma^{\beta}\rangle\hspace{-2pt}\rangle_{1}
+\frac{89}{210}\langle\hspace{-2pt}\langle\gamma^{\alpha}\rangle\hspace{-2pt}\rangle_{1}\langle\hspace{-2pt}\langle\gamma_{\alpha}\mathcal{W}\gamma^{\beta}\gamma^{\mu}\rangle\hspace{-2pt}\rangle_{0}\langle\hspace{-2pt}\langle\gamma_{\beta}\rangle\hspace{-2pt}\rangle_{1}\langle\hspace{-2pt}\langle\gamma_{\mu}\rangle\hspace{-2pt}\rangle_{1}
\nonumber
\\&-\frac{1}{210}\langle\hspace{-2pt}\langle\gamma^{\alpha}\rangle\hspace{-2pt}\rangle_{1}\langle\hspace{-2pt}\langle\gamma_{\alpha}\gamma^{\beta}\{\gamma_{\beta}\circ \mathcal{W}\}\rangle\hspace{-2pt}\rangle_{1}+\frac{1}{140}\langle\hspace{-2pt}\langle \mathcal{W}\gamma^{\alpha}\gamma^{\beta}\rangle\hspace{-2pt}\rangle_{1}\langle\hspace{-2pt}\langle\{\gamma_{\alpha}\circ\gamma_{\beta}\}\rangle\hspace{-2pt}\rangle_{1}
\nonumber
\\&+\frac{23}{140}\langle\hspace{-2pt}\langle\gamma^{\alpha}\gamma^{\beta}\rangle\hspace{-2pt}\rangle_{1}\langle\hspace{-2pt}\langle\gamma_{\alpha}\gamma_{\beta}\mathcal{W}\gamma^{\mu}\rangle\hspace{-2pt}\rangle_{0}\langle\hspace{-2pt}\langle\gamma_{\mu}\rangle\hspace{-2pt}\rangle_{1}
-\frac{3}{140}\langle\hspace{-2pt}\langle\gamma^{\alpha}\gamma^{\beta}\rangle\hspace{-2pt}\rangle_{1}\langle\hspace{-2pt}\langle\{\gamma_{\alpha}\circ\gamma_{\beta}\}\mathcal{W}\rangle\hspace{-2pt}\rangle_{1}
\nonumber
\\& - \frac{1}{{4480}} \langle\hspace{-2pt}\langle \mathcal{W}{\gamma ^\alpha } \rangle\hspace{-2pt}\rangle_{1} \langle\hspace{-2pt}\langle {\gamma _\alpha }{\gamma _\beta }{\gamma ^\beta }{\gamma _\mu }{\gamma ^\mu } \rangle\hspace{-2pt}\rangle_{0} + \frac{{13}}{{8064}} \langle\hspace{-2pt}\langle {\gamma ^\alpha } \rangle\hspace{-2pt}\rangle_{1} \langle\hspace{-2pt}\langle {\gamma _\alpha }\mathcal{W}{\gamma _\beta }{\gamma ^\beta }{\gamma _\mu }{\gamma ^\mu } \rangle\hspace{-2pt} \rangle_{0}
\nonumber
\\& - \frac{1}{{2240}} \langle\hspace{-2pt}\langle \mathcal{W}{\gamma ^\alpha }{\gamma ^\beta } \rangle\hspace{-2pt}\rangle_{1} \langle\hspace{-2pt}\langle {\gamma _\alpha }{\gamma _\beta }{\gamma _\mu }{\gamma ^\mu } \rangle\hspace{-2pt}\rangle_{0} + \frac{{41}}{{6720}} \langle\hspace{-2pt}\langle {\gamma ^\alpha }{\gamma ^\beta } \rangle\hspace{-2pt}\rangle_{1} \langle\hspace{-2pt}\langle {\gamma _\alpha }{\gamma _\beta }\mathcal{W}{\gamma _\mu }{\gamma ^\mu } \rangle\hspace{-2pt}\rangle_{0}
\nonumber
\\& + \frac{1}{{53760}} \langle\hspace{-2pt}\langle \mathcal{W}{\gamma ^\alpha }{\gamma _\alpha }{\gamma _\beta }{\gamma ^\beta }{\gamma _\mu }{\gamma ^\mu } \rangle\hspace{-2pt}\rangle_{0} - \frac{1}{{210}} \langle\hspace{-2pt}\langle \{ \mathcal{W} \circ {\gamma ^\alpha }\}  \rangle\hspace{-2pt}\rangle_{1} \langle\hspace{-2pt}\langle {\gamma _\alpha }{\gamma _\beta }{\gamma ^\beta } \rangle\hspace{-2pt}\rangle_{1}
\nonumber
\\& - \frac{1}{{5760}} \langle\hspace{-2pt}\langle \mathcal{W}{\gamma ^\alpha }{\gamma _\alpha }\{ {\gamma _\beta } \circ {\gamma ^\beta }\}  \rangle\hspace{-2pt}\rangle_{1} - \frac{1}{{2688}} \langle\hspace{-2pt}\langle {\gamma ^\alpha }{\gamma _\alpha }{\gamma ^\beta } \rangle\hspace{-2pt}\rangle_{1} \langle\hspace{-2pt}\langle {\gamma _\beta }\mathcal{W}{\gamma _\mu }{\gamma ^\mu } \rangle\hspace{-2pt}\rangle_{0}
\nonumber
\\& - \frac{1}{{5040}} \langle\hspace{-2pt}\langle {\gamma ^\alpha }{\gamma _\alpha }{\gamma ^\beta }\{ {\gamma _\beta } \circ \mathcal{W}\}  \rangle\hspace{-2pt}\rangle_{1} + \frac{1}{{3780}} \langle\hspace{-2pt}\langle \mathcal{W}{\gamma _\alpha }{\gamma _\beta }{\gamma _\mu } \rangle\hspace{-2pt}\rangle_{1} \langle\hspace{-2pt}\langle {\gamma ^\alpha }{\gamma ^\beta }{\gamma ^\mu } \rangle\hspace{-2pt}\rangle_{0}
\nonumber
\\& + \frac{1}{{252}} \langle\hspace{-2pt}\langle {\gamma _\alpha }{\gamma _\beta }{\gamma _\mu } \rangle\hspace{-2pt}\rangle_{1} \langle\hspace{-2pt}\langle \mathcal{W}{\gamma ^\alpha }{\gamma ^\beta }{\gamma ^\mu } \rangle\hspace{-2pt}\rangle_{0}.
\label{genus3TRR}
\end{align}
By equations~\eqref{genus3TRR}, \eqref{genus0TRR}, \eqref{genus1TRR}, \eqref{genus2TRR} and their derivatives, both sides of equation~\eqref{g3LP} can be reduced to the following expression
 \begin{align*}
&\frac{7}{5760}\langle\hspace{-2pt}\langle \{\mathcal{W}\circ\Delta\circ\Delta\}\rangle\hspace{-2pt}\rangle_{1}
 +\frac{11}{2903040}\langle\hspace{-2pt}\langle \{\mathcal{W}\Delta\Delta\Delta\}\rangle\hspace{-2pt}\rangle_{1}
  \\&+\frac{19}{967680}\langle\hspace{-2pt}\langle \{\mathcal{W}\circ\Delta\}\Delta\gamma^{\alpha}\gamma_{\alpha}\rangle\hspace{-2pt}\rangle_{0}
  +\frac{1}{120960}\langle\hspace{-2pt}\langle \mathcal{W}\{\Delta\circ\Delta\}\gamma^{\alpha}\gamma_{\alpha}\rangle\hspace{-2pt}\rangle_{0}
  \\&+\frac{1}{60480}\langle\hspace{-2pt}\langle \{\mathcal{W}\circ\gamma^{\alpha}\}\gamma_{\alpha}\Delta\Delta\rangle\hspace{-2pt}\rangle_{0}
   +\frac{1}{11520}\langle\hspace{-2pt}\langle \{\mathcal{W}\circ\Delta\circ\gamma^{\alpha}\}\gamma_{\alpha}\gamma^{\beta}\gamma_{\beta}\rangle\hspace{-2pt}\rangle_{0}.
 \end{align*}
 \subsection{Case of genus $g=4$}
  Due to the dimension constraints, we only need to equation \eqref{eqXPrel} holds for $r=0$, $r=1$ and $r=2$.  For simplicity, we only prove the case for $r=0$ and the other case can be similarly done, i.e.
    \begin{align}
& \langle\hspace{-2pt}\langle T^{8}(\mathcal{W})\rangle\hspace{-2pt}\rangle_{4}\nonumber
   \\=& - \frac{2}{4}{\langle\hspace{-2pt}\langle \mathcal{W}{T^5}({\gamma ^\alpha })\rangle\hspace{-2pt}\rangle _2}{\langle\hspace{-2pt}\langle {T^2}({\gamma _\alpha })\rangle\hspace{-2pt}\rangle _2} + \frac{1}{4}{\langle\hspace{-2pt}\langle \mathcal{W}{T^6}({\gamma ^\alpha })\rangle\hspace{-2pt}\rangle _3}{\langle\hspace{-2pt}\langle T({\gamma _\alpha })\rangle\hspace{-2pt}\rangle _1}
   - \frac{1}{4}{\langle\hspace{-2pt}\langle \mathcal{W}{T^7}({\gamma ^\alpha })\rangle\hspace{-2pt}\rangle _3}{\langle\hspace{-2pt}\langle {\gamma _\alpha }\rangle\hspace{-2pt}\rangle _1}.
   \label{g4LP}
 \end{align}
By the equations~\eqref{g4LP}, \eqref{genus0TRR}, \eqref{genus1TRR}, \eqref{genus3TRR} and their derivatives, we can reduce both sides of equation~\eqref{g4LP} to an expression only involving genus-0 and genus-1 invariants. As the expression itself is long and not so important, we omit it here.

%
%

\vspace{40pt}

\appendix
\centerline {\bf \Large Appendix}
\section{Expression of $A(\mathcal{W})$}
\label{sec:A(W)}
In this appendix, we give the explicit expression of $A(\mathcal{W})$ in equation \eqref{genus4TRR}.


\vspace{30pt} \noindent
Xin Wang \\
School of Mathematics  \\
Shandong University, Jinan, China \\
E-mail address:{\it wangxin2015@sdu.edu.cn}

\end{document}